\documentclass[12pt]{article}
\usepackage{latexsym,amssymb,psfig}

\font\of=msbm10 scaled 1200
\def\R{\mbox{\of R}}
\def\C{\mbox{\of C}}
\def\Z{\mbox{\of Z}}

\newtheorem{definition}{Definition}

\newtheorem{theorem}{Theorem}
\newtheorem{lemma}{Lemma}
\newtheorem{proposition}{Proposition}

\title{ Complete hyperelliptic integrals of the first kind and their
non-oscillation }

\author{Lubomir Gavrilov \\
\normalsize \it Laboratoire Emile Picard , CNRS UMR 5580,  Universit\'e
Paul Sabatier\\
\normalsize \it 118, route de Narbonne, 31062 Toulouse Cedex, France \\
\\ Iliya D. Iliev\\
\normalsize \it Institute of Mathematics, Bulgarian Academy of Sciences\\
\normalsize \it P.O. Box
373, 1090 Sofia, Bulgaria } \date{}

\begin{document}
\maketitle
\begin{abstract}
Let $P(x)$ be a real polynomial of degree $2g+1$, $H=y^2+P(x)$ and
$\delta(h)$ be an oval contained in the level set $\{H=h\}$.
We study complete Abelian integrals of the form
$$I(h)=\int_{\delta(h)}
\frac{(\alpha_0+\alpha_1 x+\ldots + \alpha_{g-1}x^{g-1})dx}{y},
\;\;h\in \Sigma,$$
where $\alpha_i$ are real and $\Sigma\subset \R$ is a maximal open interval
on which a continuous family of ovals $\{\delta(h)\}$ exists. We show that
the $g$-dimensional real vector space of these integrals is not Chebyshev
in general: for any $g>1$, there are hyperelliptic Hamiltonians $H$ and
continuous families of ovals $\delta(h)\subset\{H=h\}$, $h\in\Sigma$,
such that the Abelian integral $I(h)$ can
have at least $[\frac32g]-1$ zeros in $\Sigma$. Our main result is Theorem
\ref{main} in which we show that when $g=2$, exceptional families of ovals
$\{\delta(h)\}$ exist, such that the corresponding vector space is still
Chebyshev.
\end{abstract}

\section{Introduction}

Take real polynomials  $H,f,g\in \R[x,y]$ and let $\delta (h) \subset
\{(x,y)\in \R^2: H(x,y)=h\}$, $h\in\Sigma$,
be a continuous family of ovals.
For  sufficiently small $\varepsilon$ and generic $H, f, g,$
the limit cycles of the perturbed plane Hamiltonian system
$$
dH + \varepsilon (f dx + g dy) = 0
$$
which tend to certain ovals from the continuous
family when $\varepsilon\to 0$, are in one-to-one
correspondence with the zeros of the complete Abelian integral
$$
I(h) = \int_{\delta (h)} f(x,y) dx + g(x,y) dy, \quad h\in\Sigma.
$$
For this reason the problem of finding the zeros of $I(h)$ in terms of the
degrees of $H,f,g$
was called  by Arnold \cite[p. 313]{Arnold1} the
``{\it weakened 16th Hilbert problem}"  (compare to Hilbert \cite{hil}, see
also Arnold \cite{Arnold2,Arnold1+,ten}).
Note that the level sets $ \{H=h\} $ will contain in general several
continuous families of ovals which need be considered separately.

It follows from the Varchenko-Khovanskii theorem that the number of the
zeros of $I(h)$ is bounded
by a constant $N(n,d)<\infty$, uniformly in all $H,f,g$, such that
$deg\,f\leq d$,  $deg\,g \leq d$,  $deg\,H\leq n$.

In the so called ``hyperelliptic case" ($H=y^2+P(x)$) it was proved by
Novikov and Yakovenko\cite{nov-yak} that there exists an algorithm producing
a function $C(n,d)$ such that $\infty > C(n,d)> N(n,d)$,
which is given by a tower function (an iterated exponent) of
height at least five. Their
proof is based on the analytic properties of a suitable Picard-Fuchs system
satisfied by Abelian
integrals (including the magnitude of the coefficients of the system and some
restrictions on its monodromy group).

The progress in solving the weakened 16th Hilbert problem (finding the
{\it exact} number of the zeros of
$I(h)$) concerned so far the ``elliptic case" only (the complex algebraic
curve $\{H=h\}$
is of genus at most one). It was proved for instance
 that in several cases,
the vector space ${\cal A}_{H,d}$
of Abelian integrals of degree $d$ polynomials along the ovals of $H$,
obeys the so called {\it Chebyshev property} (the number
of the zeros of each integral is smaller than the dimension of the vector
space ${\cal A}_{H,d}$), see \cite{Petrov86,gav98,gav99,gas-li}.
In this relation Arnold asked in \cite[the 7th problem]{ten} whether the
$g$-dimensional vector
space of Abelian integrals
\begin{equation}
\label{question}
I(h)=\int_{\delta(h)}
\frac{(\alpha_0+\alpha_1 x+\ldots + \alpha_{g-1}x^{g-1})dx}{y},
\;\;h\in \Sigma
\end{equation}
 where
$
\delta (h)\subset \{(x,y)\in \R^2:\; y^2+P(x)=h\},\;\;g=[\frac12(deg\, P-1)]>1,
$
is Chebyshev.

In an attempt to solve this problem Givental \cite{giv} obtained a
non-oscillation theorem for Lagrangian
planes of the Picard-Fuchs system satisfied by the Abelian integrals
$\int_{\delta(h)}x^idx/y$. He
used the fact that every Picard-Fuchs system has a Hamiltonian (or Poisson)
structure. He failed,
however, to produce bounds for the zeros of the Abelian integrals.

In the present paper we try to explore the algebro-geometrical properties
of the Abelian integrals,
without using the Picard-Fuchs system. Our main result is Theorem
\ref{main} in which we show that
when $g=2$ and $deg\,P=5$, exceptional families of ovals
$\{\delta(h)\}$, $ h\in \Sigma$, exist such that every Abelian integral of
the form
$$
I(h)=\int_{\delta(h)}
\frac{(\alpha_0+\alpha_1 x)dx}{y},\quad \alpha _0^2+\alpha _1^2\neq 0
$$
has at most one simple zero on the  interval $\Sigma$. At a first sight
this seems to be an easy
observation (equivalent to the claim that $(\int_{\delta(h)}
xdx/y)/(\int_{\delta(h)}dx/y)$ is
monotonous on $\Sigma$). This is, however, the first result of such a kind
for non-elliptic curves
$\{H=h\}$. The proof uses the Riemann bilinear relations on differentials
of the first kind together
with  the fact that a Jacobian variety with its polarization cannot be a
direct product of
principally polarized Abelian varieties. Our arguments can be adapted to
other situations
(for instance when the degree of $P(x)$ is $6$), but we shall not do
this here.

We give also a negative answer (Proposition \ref{nocheb}) to the initial
question posed by Arnold:
it turns out that there exist Abelian integrals of the form
(\ref{question}) with exactly
$[\frac32g]-1$ zeros in a neighborhood of the origin.

\section{The cyclicity at a center}
In this section we prove that in general, the integral (\ref{question})
does not belong to a Chebyshev space. Consider first the particular case
when $n=5$. Let $H= y^2+ P(x)$ where $P(x)= x^2(1+a_1 x +a_2 x^2+a_3 x^3)$,
$a_1,a_2, a_3 \in \R$, $a_3\neq 0$, and denote by
$\delta (t) \subset \{y^2+P(x)=t\}$  a continuous family of cycles vanishing
at the origin as $t$ tends to zero.

\begin{proposition}
\label{cyclicity}
The Abelian integral
$$
I(t) = \int_{\delta (t)} \frac{(\alpha_0+\alpha_1 x) dx}{y}=
\alpha_0 I_0(t)+\alpha_1 I_1(t)\not\equiv 0
$$
can have a zero at the origin of multiplicity at most two. Moreover, for
fixed $\alpha_1, a_2, a_3$, $\alpha_1 a_3\neq 0$, there exist a sufficiently
small
$\varepsilon >0$  and $\alpha_0$, $a_1$ in a small neighborhood of zero,
such that $I(t)$ has exactly two simple zeros in the interval
$(0,\varepsilon )$.
\end{proposition}
The above shows that the real vector space generated by the functions
$I_0(t)$, $I_1(t)$ is not Chebyshev. Recall that a real vector space $V$ of
functions defined on some interval $\Sigma $ is said to be Chebyshev,
provided that each $f\in V$ has at most $\dim V - 1$ zeros (counted with
multiplicity) on $\Sigma $, and Chebyshev with accuracy $m$, if
each $f\in V$ has at most $\dim V+m - 1$ zeros there.

\vspace{2ex}
\noindent
{\bf Proof of Proposition \ref{cyclicity}.}
For a small $x$, denote $X=x(1+a_1 x +a_2 x^2+a_3 x^3)^{1/2}$.
Then an easy calculation yields the inverse transformation
$$\textstyle
x=\varphi(X)\sim X-\frac12a_1X^2+(\frac58a_1^2-\frac12a_2)X^3
-(a_1^3-\frac32a_1a_2+\frac12a_3)X^4+\ldots \quad \mbox{\rm as}\quad
X\to 0.$$
Therefore, for small positive $t$,
$$\begin{array}{l}{\displaystyle
I(t)=\int_{y^2+X^2=t}\frac{[\alpha_0+\alpha_1\varphi(X)]
\varphi'(X)}{y}dX}\\[5mm]
{\displaystyle
= \int_{y^2+X^2=t} \frac{\alpha_0[1+O(X^2)] - \alpha_1[\frac32a_1X^2
+\frac52(\frac{21}{8}a_1^3-\frac72a_1a_2+a_3)X^4 +O(X^6)]}{y}dX}\\[5mm]
=2\pi\alpha_0[1+O(t)]-2\pi\alpha_1[\frac34a_1t+
\frac{15}{16}(\frac{21}{8}a_1^3- \frac72a_1a_2+a_3)t^2 +O(t^3)].
\end{array}$$
When $\alpha_0=a_1=0$, the function $I(t)$ has a double zero at the
origin since $\alpha_1 a_3\neq 0$.
Taking $\alpha_1=1$, $|\alpha_0|<\!\!<|a_1|<\!\!<|a_3|$ and
$\alpha_0 a_1>0>a_1a_3$, one finishes the proof of the proposition. $\Box$

Clearly, the above approach could be applied to hyperelliptic Hamiltonians
of any degree that have a center.

Denote by ${\cal H}_n$ the set of Hamiltonians
$H=y^2+x^2+a_1 x^3+\ldots +a_{n-2}x^n$, where $a_k\in \R$, $a_{n-2}\neq 0$.
Denote by $g=[\frac12(n-1)]\geq 2$ the genus of the hyperelliptic curve.
Given $H\in {\cal H}_n$, let $\delta(t)$, $t\in (0, T_H)$, be an oval from
the continuous family of ovals surrounding the center at the origin.
Consider the vector space ${\cal A}_H$ formed by the Abelian integrals

\begin{equation}
\label{g}
I(t) = \int_{\delta (t)} \frac{G(x) dx}{y},\quad t\in(0,T_H),
\quad G\in \R[x],\quad deg\, G<g.
\end{equation}
Let us point out that $dim\,{\cal H}_H=g$ in general, but
$dim\,{\cal H}_H=[\frac12(g+1)]$ provided that $H$ contains no odd degree
monomials. In the proposition below, we give a partial answer to the question
raised by V. Arnold'd in \cite{ten}.

\begin{proposition}
\label{nocheb}
There exist $H\in{\cal H}_n$ and a polynomial $G$ such that the
corresponding integral $(\ref{g})$ has exactly $[\frac32g]-1$ small positive
simple zeroes. Thus, such a space ${\cal A}_H$ could be Chebyshev with
accuracy at least $[\frac12g]$ in $(0,T_H)$.
\end{proposition}

Given an integer $k\geq 0$ and $H\in{\cal H}_n$, denote
$$I_k(t)=\int_{\delta(t)}\frac{x^kdx}{y}, \quad 0<t<T_H.$$
The proof of the above proposition is based on the following lemma.

\begin{lemma}
\label{int} The following asymptotic expansions hold for
small positive $t$:

\vspace{1ex}
\noindent
{\em(i)}$\quad I_{2k}(t)= t^k[c_k+O(t)],\;\;\mbox{\it where} \;\;
c_k=2\pi\frac{(2k-1)!!}{(2k)!!},$

\vspace{1ex}
\noindent
{\em(ii)}$\quad\displaystyle I_{2k+1}(t)=-{\textstyle\frac12}
\sum\limits_{j=1}^g a_{2j-1}(2k+2j+1)t^{k+j}[c_{k+j}+O(t)],$

\vspace{1ex}
\noindent
therefore the coefficients in the last expansion belong to the ideal
generated by the odd-numbered parameters $\{a_1, a_3,\ldots, a_{2g-1}\}$.
\end{lemma}

\noindent
{\bf Proof.} Take as above $X=x(1+a_1 x +\ldots+ a_{n-2} x^{n-2})^{1/2}$,
then the inverse transformation $x=\varphi(X,a_1,a_2,\ldots,a_{n-2})$
for small $X$ reads $ x=X+A_1X^2+A_2X^3 +A_3X^4+\ldots$,
where the coefficients $A_j=A_j(a_1,a_2,\ldots,a_{n-2})$ are polynomials.
It is easy to see that if one attaches a weight $k$ to the coefficient $a_k$,
then $A_j$ is a polynomial of weight $j$. Indeed, directly from the definition
of $\varphi$ it follows that
$$\varphi(\omega X,a_1,a_2,\ldots,a_{n-2})=\omega\varphi(X, \omega a_1,
\omega^2 a_2,\ldots, \omega^{n-2}a_{n-2}).$$
Applying this identity to
$$\varphi(X, a_1, \ldots, a_{n-2})=X(1+\sum_{k=1}^\infty A_k(a_1,\ldots,
a_{n-2})X^k),$$
we get $A_k(\omega a_1, \ldots, \omega^{n-2} a_{n-2})=\omega^k
A_k(a_1, \ldots, a_{n-2})$.
To calculate the coefficient at $a_k$ in $A_k$, we take $a_1=\ldots=a_{k-1}=0$,
and easily obtain that it equals $-\frac12$. Therefore, one deduces that
$A_k=-\frac12a_k+A_k^*$ where $A_k^*$ belongs to the polynomial ideal
generated by $a_1,\ldots, a_{k-1}$.

In case (i), one immediately obtains that for small positive $t$,
$$I_{2k}(t)=\frac{1}{2k+1}\int_{y^2+X^2=t}\frac{d\varphi^{2k+1}(X)}{y}=
\int_{y^2+X^2=t}\frac{(X^{2k}+\ldots)dX}{y}=c_kt^k+O(t^{k+1}).
$$
To calculate the integral in case (ii), we write $\varphi^{2k+2}(X)$
in the form
$$\varphi^{2k+2}=X^{2k+2}\left(1+\sum\limits_{j=1}^\infty
(-{\textstyle\frac12}a_j+A_j^*)X^j\right)^{2k+2}$$
$$=X^{2k+2}\left(1+\sum\limits_{j=1}^\infty
[(2k+2)(-{\textstyle\frac12}a_j+A_j^*)X^j+O(X^{j+1})]\right)$$
$$=X^{2k+2}\left(1-\sum\limits_{j=1}^{n-2} (k+1)a_jX^j[1+O(X)]\right)$$
$$=-(k+1)\sum\limits_{j=1}^g a_{2j-1}X^{2k+2j+1}[1+O(X^2)]+ \Phi(X)$$
where $\Phi(X)$ is a series in $X$ containing only even order powers.
Therefore,
$$I_{2k+1}(t)=\frac{1}{2k+2}\int_{y^2+X^2=t}\frac{d\varphi^{2k+2}(X)}{y}$$
$$=-{\textstyle\frac12}\sum\limits_{j=1}^g a_{2j-1}(2k+2j+1)t^{k+j}
[c_{k+j}+O(t)].\;\;\Box$$

\vspace{2ex}
\noindent
{\bf Proof of Proposition \ref{nocheb}.} Denote $m=[\frac12g]$ and take
in (\ref{g})
$$G(x)=\sum_{k=0}^{m-1}\gamma_{2k}x^{2k}-x^{2m-1}.$$
Then by Lemma \ref{int},
$$I(t)=\sum_{k=0}^{m-1}\gamma_{2k}t^k[c_k+O(t)]+{\textstyle\frac12}
\sum_{j=1}^g a_{2j-1}(2m+2j-1)t^{m+j-1}[c_{m+j-1}+O(t)].$$
Now, freeze the even-numbered coefficients in $H$ and then choose $\gamma_{2j}$
and $a_{2j-1}$ so that $\gamma_{2j}\gamma_{2j+2}<0$,
$a_{2j-1}a_{2j+1}<0$, $\gamma_{2m-2}a_1<0$ and
$$|\gamma_0|<\!\!<|\gamma_2|<\!\!<\ldots<\!\!<|\gamma_{2m-2}|
<\!\!<|a_1|<\!\!<\ldots <\!\!<|a_{2g-1}|<\!\!<1.$$
Such a choice guarantees that $I(t)$ will have exactly $m+g-1=[\frac32g]-1$
zeroes in a certain small interval $(0,\varepsilon)$. $\Box$

This result raises the problem to describe the hyperelliptic Hamiltonians
$H$ and the continuous families of ovals $\{\delta(h)\}$ for which the
space of integrals (\ref{question}) obeys the Chebyshev property.
In the rest of the paper, we concentrate our efforts on the simplest case
when $deg\,P=5$. We begin with some preliminaries.

\section{The normal form, the bifurcation diagram and the Dynkin diagram
for $n=5$}
We consider the hyperelliptic Hamiltonian $H=\frac12y^2+P(x)$
where $P\in \R[x]$ and $deg\,P=5$. Assume that a certain level set
$\{H=h\}$ contains an oval (compact closed smooth curve without critical
points). As the Poincar\'e index of the respective Hamiltonian vector
field $X_H$ is zero, then $X_H$ has at least two different real critical
points. Without any loss of generality, we can place all the four critical
points at $(0,0)$, $(\mu, 0)$, $(\lambda,0)$ and $(1,0)$, where either
$0\leq\mu\leq\lambda\leq 1$ (the real case)
or $\bar\mu=\lambda\in\C\setminus\R$ (the complex case).
Finally, we can assume that the coefficient at $x^5$ is positive.
Then after rescaling of the $y$ variable, we come to the normal form
\begin{equation}
H(x,y)=\frac12y^2-\frac{\lambda\mu}{2}x^2+
\frac{\lambda+\mu+\lambda\mu}{3}x^3-\frac{1+\lambda+\mu}{4}x^4+\frac15x^5
\equiv \frac12y^2+P(x)
\label{normal}
\end{equation}
and the Hamiltonian flow $X_H$ takes the form
\begin{equation}
\label{system}
\begin{array}{l}
\dot{x}=y,\\
\dot{y}=-x(x-\mu)(x-\lambda)(x-1).
\end{array}
\end{equation}
Until the end of the paper, we will use this normal form of $H$ only.

\begin{figure}[htbp]
\begin{center}
  \psfig{file=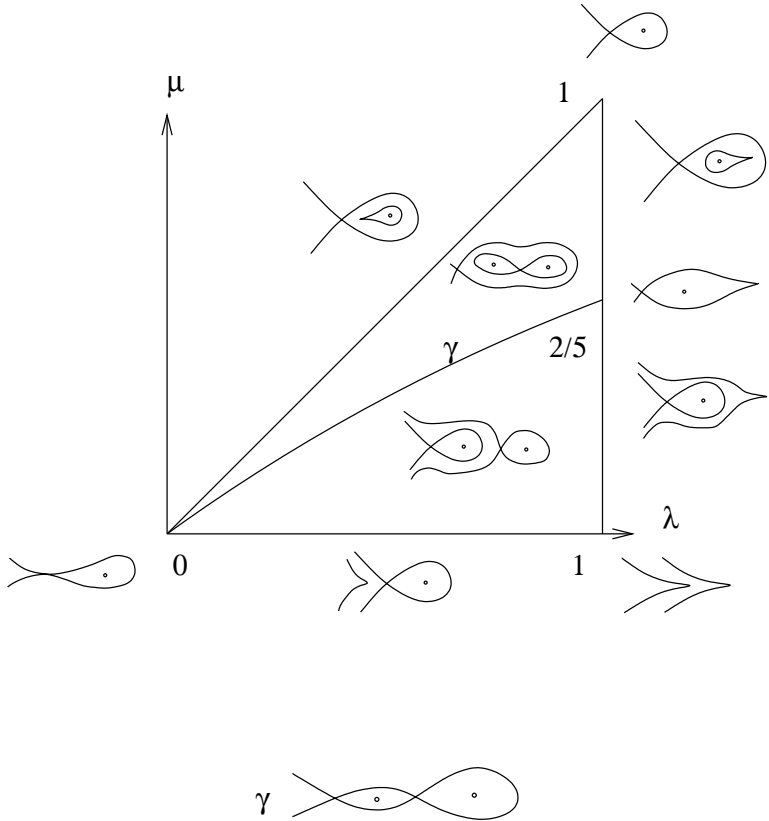}
\end{center}
\caption{ Bifurcation diagram in the $(\lambda,\mu)$-plane in the real case}
\label{bdrc}
\end{figure}

Clearly, the origin is a hyperbolic saddle if $\mu\neq 0$, and $(1,0)$ is
a nondegenerate center provided that $\lambda\neq 1$. In addition,
when $0<\mu<\lambda<1$ in the real case, $(\mu, 0)$ is a center
and $(\lambda, 0)$ is a saddle. Denote by $h_0, h_\mu, h_\lambda$ and $h_1$
the corresponding critical levels of $H$. One has
$$\begin{array}{l}
h_0=0,\quad h_1=-\frac{1}{60}(3-5\lambda-5\mu+10\lambda\mu),\\[2mm]
h_\lambda=-\frac{\lambda^3}{60}(3\lambda^2-5\lambda\mu-5\lambda+10\mu),\\[2mm]
h_\mu=-\frac{\mu^3}{60}(3\mu^2-5\lambda\mu-5\mu+10\lambda). \end{array}$$
If $0<\mu<\lambda<1$, we have $h_0>h_\mu$, $h_\lambda>h_\mu$, $h_\lambda>h_1$.
The bifurcation diagram $D$ in the real case consists of the boundary
of the triangle $T=\{0\leq\mu\leq\lambda\leq 1\}$ and the curve $\gamma:\;
h_\lambda=h_0$ having an equation
$$\mu=\frac{3\lambda^2-5\lambda}{5(\lambda-2)},$$
see Fig. \ref{bdrc}.

There are two open components in $T\setminus D$ giving the
two generic cases. Namely, for parameters $\lambda, \mu$ above the
curve $\gamma$, $X_H$ has three period annuli, and below this curve $X_H$
has just two period annuli. There are also five codimension-one and four
codimension-two degenerate (nongeneric) cases corresponding
to the parameters on the bifurcation diagram $D$. In the complex case,
there is a unique period annulus around $(1,0)$ which terminates
at the saddle-loop through $(0,0)$ and the phase portraits for all
$\lambda\in\C\setminus\R$ are topologically equivalent. We will denote
the period annuli (continuous families of ovals) around $(1,0)$, $(\mu,0)$
and the eight loop by ${\cal O}_1$, ${\cal O}_\mu$ and ${\cal O}_e$,
respectively.

Below we describe the Dynkin diagram of the real polynomial $H(x,y)$
when the critical levels are distinct. Recall that it is a graph
defined in the same way as the Dynkin diagram (or D-diagram)
of a germ of an analytical function with an isolated singularity, see
\cite{Arnold,dim}. The vertices of the diagram are in one-to-one
correspondence with the vanishing cycles of the polynomial. Two
vertices are connected by an edge if the intersection number of the
cycles is not zero. To an edge one associates also a label equal
to the intersection number, but it will be omitted here because is of no
importance for us. The vanishing cycles ``vanish" along suitable paths,
and hence the Dynkin diagram depends also on the family of paths. The
precise definition in our situation is as follows.

\begin{figure}
\setlength{\unitlength}{.03\textwidth}
\begin{picture}(30,15)(0,0)
\put(0,5){
 \begin{picture}(15,15)(0,0)
 \put(0,0){\circle*{.3}} \put(0,0.5){$h_1$}
 \put(3,0){\circle*{.3}}  \put(3,0.5){$h_\mu$}
 \put(6,0){\circle*{.3}} \put(6,0.5){$h_\lambda$}
 \put(9,0){\circle*{.3}} \put(9,0.5){$h_0$}
 \put(6,0){\thicklines{\vector(1,0){6}}}

 \qbezier[20](0,0)(0,5)(6,6)
 \qbezier[20](3,0)(3,2)(6,6)
 \qbezier[20](6,0)(5,3)(6,6)
 \qbezier[20](9,0)(7.5,2.5)(6,6)
 \put(6,6){\circle*{.3}} \put(6.5,6){$\tilde{h}$}
 \end{picture}}
\put(20,5){
 \begin{picture}(15,15)(0,0)
 \put(0,0){\line(1,0){3}} \put(0,0){\circle*{.3}}
 \put(3,0){\line(1,0){3}} \put(3,0){\circle*{.3}}
 \put(6,0){\line(1,0){3}} \put(6,0){\circle*{.3}}
 \put(9,0){\circle*{.3}}

 \put(0,.5){$\delta _1$}
 \put(2.5,.5){$\delta _\lambda$}
 \put(5.5,.5){$\delta _\mu$}
 \put(8.5,.5){$\delta _0$}

 \end{picture}
              }
\put(7,0){ (i)}
\put(22,0){ (ii)}
\end{picture}
\caption{
(i) The paths $l_0$, $l_1$, $l_\lambda$, $l_\mu$ in the real case
$0<\mu<\lambda<1$ when $h_1<h_\mu<h_\lambda<h_0$;
(ii) Dynkin diagram in the real case $0<\mu <\lambda <1$ }
\label{realdynk}
\end{figure}
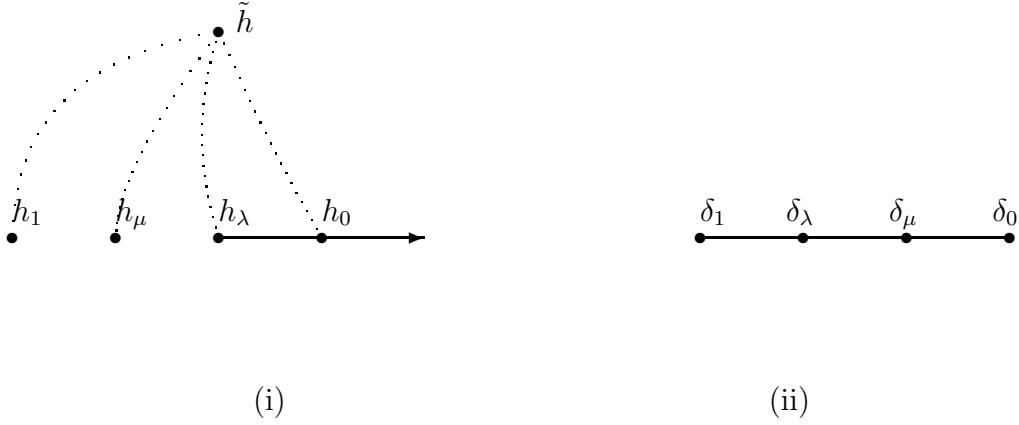

Consider first the real case $0<\mu<\lambda<1$. There are five possible
distributions of the critical values as follows:
$$\begin{array}{lll}
h_1 < h_\mu < h_\lambda < h_0, &
h_1 < h_\mu < h_0 < h_\lambda, &
h_\mu < h_1 < h_0 < h_\lambda, \\
h_\mu < h_0 < h_1 < h_\lambda, &
h_\mu < h_1 < h_\lambda < h_0. &
\end{array}$$
For definiteness, assume that $h_1<h_\mu < h_\lambda < h_0$.
Denote ${\cal D} = \C \backslash [h_\lambda , \infty)$ and let
$$
l_0,l_1,l_\lambda,l_\mu:\; [0,1] \rightarrow \{h\in \C: Im\,h \geq 0\}
$$
be continuous paths connecting some fixed regular value
$\tilde{h} \in {\cal D}$ to
$h_0$, $h_1$, $h_\lambda$, $h_\mu$ respectively, as shown on
Fig. \ref{realdynk}.

%

The Dynkin diagram of $H$ with respect to the above paths
is a graph with four vertices corresponding to the
four families of cycles vanishing along the paths
$l_0$, $l_1$, $l_\lambda$, $l_\mu$, which we denote by
$\delta_0(h)$, $\delta_1(h)$, $\delta_\mu (h)$, $\delta_\lambda (h)$
respectively. Two vertices are connected by an
edge if the intersection number of the cycles for $h=\tilde{h}$
is not zero (in which case it is $\pm 1$ and may be supposed equal
to one). The computation of the Dynkin diagram of a real polynomial
with real critical points only is well known, see for instance
\cite{cam,Arnold}. It yields the diagram shown on Fig. \ref{realdynk}.
We note that the computation depends only on the distribution of the
four critical points $0<\mu<\lambda<1$ and not on the distribution of
the respective critical values (hence, the Dynkin diagram is the same
for all possible arrangements of the critical values
with corresponding paths $l_0$, $l_1$, $l_\lambda$, $l_\mu$).

Consider now the complex case when $\lambda = \bar{\mu}\not\in\R $.
The method we use is a deformation of $H$ to a polynomial having real
critical points. Denote ${\cal D}=\C\backslash [h_0,\infty)$ and let
$$
l_0, l_1, l_\lambda, l_\mu:\; [0,1] \rightarrow {\cal D}
$$
be continuous paths connecting some fixed regular value
$\tilde{h} \in {\cal D}$ to $h_0$, $h_1$, $h_\lambda$, $h_\mu$
respectively, as shown on Fig. \ref{dynkin1}. The Dynkin diagram
related to these paths is defined as above. Consider in
the complex plane the hyperbola $\Gamma$ defined by the equation
\begin{equation}
\label{hyperbola}
\Gamma:\;\;
h_\lambda=h_\mu\quad \Leftrightarrow\quad (Im\,\lambda)^2=
5Re\,\lambda(Re\,\lambda-1)
\end{equation}
and drawn on Fig. \ref{imconj}.

Suppose for definiteness that $\lambda$ lies inside its left branch.
First, by continuous deformation of $\lambda = \bar{\mu }$ such that
in the course of the deformation $\lambda $ remains in the left branch
of (\ref{hyperbola}), we may achieve that $\lambda =\mu<0<1$.
This, combined with $h_\lambda =h_\mu <0$ already implies that
\begin{itemize}
\item
the intersection number of $\delta _1(\tilde{h})$ and $\delta _0(\tilde{h})$
is equal to $\pm 1$;
\item
the intersection number of $\delta _\lambda (\tilde{h})$,
$\delta _\mu  (\tilde{h})$ with $\delta _1 (\tilde{h})$ is zero;
\item
the intersection number of $\delta _\lambda (\tilde{h})$ and
$\delta _\mu  (\tilde{h})$ is equal to $\pm 1$;
\item
The intersection number of either $\delta _\lambda (\tilde{h})$ or
$\delta _\mu(\tilde{h})$ with $\delta _0(\tilde{h})$ is equal
to $\pm 1$.
\end{itemize}
Finally, we may use the real structure of the complex algebraic  curve
$$\Gamma _h = \{(x,y)\in \C^2 :\; H(x,y) =h\}$$
for $h\in \R$.
This structure is defined by an antiholomorphic involution $i$ on the
curve which is the complex conjugation. It induces on its hand an
involution $i_*$ of $H_1(\Gamma _h, \Z)$, $(i_* )^2=id$.
It is easy to see that for $h<0$ we have
$$
i_* \delta _\lambda(h) = \delta _\mu(h) ,\quad
i_* \delta _0(h) = - \delta _0(h) .
$$
This combined with the identity
$$
\overline{\langle \delta (h) , \gamma (h)\rangle} =
- \langle i_*\delta (h) , i_*\gamma (h)\rangle
$$
shows that
$$\langle \delta_\lambda  (h) , \delta _0 (h) \rangle=
\langle \delta_\mu   (h) , \delta _0 (h)\rangle.
$$
This implies the Dynkin diagram as shown on Fig. \ref{dynkin1}(ii)).
The remaining cases shown there are studied in a similar way.

\begin{figure}
\setlength{\unitlength}{.03\textwidth}

\begin{picture}(30,15)(0,0)
\put(0,5){
   \begin{picture}(15,15)(0,0)
\put(0,0){\circle*{.3}} \put(0,0.5){$h_1$}
\put(3,0){\circle*{.3}}  \put(3,0.5){$h_0$}
\put(1.5,1.5){\circle*{.3}} \put(1.5,2){$h_\lambda$}
\put(1.5,-1.5){\circle*{.3}} \put(2,-1.5){$h_\mu$}
\put(3,0){\thicklines{\vector(1,0){9}}}

\qbezier[20](0,0)(0,5)(6,6)
\qbezier[20](3,0)(3,2)(6,6)
\qbezier[20](1.5,1.5)(4,5)(6,6)
\qbezier[25](1.5,-1.5)(3,3.5)(6,6)
\put(6,6){\circle*{.3}} \put(6.5,6){$\tilde{h}$}
   \end{picture}}
\put(20,5){
   \begin{picture}(15,15)(0,0)
\put(0,0){\line(1,0){3}} \put(0,0){\circle*{.3}}
\put(3,0){\line(1,1){2.1213}} \put(3,0){\circle*{.3}}
\put(3,0){\line(1,-1){2.1213}}
\put(5.1213,2.1213){\circle*{.3}}
\put(5.1213,-2.1213){\circle*{.3}}

\put(2.5,0.5){$\delta _0$}
\put(0,0.5){$\delta _1$}
\put(5.6213,2.1213){$\delta _\lambda$}
\put(5.6213,-2.1213){$\delta _\mu$}

\put(5.1213,-2.1213){\line(0,1){4.2426}}
   \end{picture}
              }
\put(7,0){ (i)}
\put(22,0){ (ii)}

\end{picture}
\begin{picture}(30,15)(0,0)
\put(3,7){
 \begin{picture}(15,15)(0,0)
 \put(0,0){\line(1,0){3}} \put(0,0){\circle*{.3}}
 \put(3,0){\line(1,1){2.1213}} \put(3,0){\circle*{.3}}
 \put(3,0){\line(1,-1){2.1213}}
 \put(5.1213,2.1213){\circle*{.3}}
 \put(5.1213,-2.1213){\circle*{.3}}

 \put(2.5,0.5){$\delta _1$}
 \put(0,0.5){$\delta _0$}
 \put(5.6213,2.1213){$\delta _\lambda$}
 \put(5.6213,-2.1213){$\delta _\mu$}

 \put(5.1213,-2.1213){\line(0,1){4.2426}}
    \end{picture}
          }
\put(20,5){
 \begin{picture}(15,15)(0,0)
 \put(0,0){\line(1,0){4}} \put(0,0){\circle*{.3}}
 \put(4,0){\line(0,1){4}} \put(4,0){\circle*{.3}}
 \put(0,4){\line(1,0){4}}
 \put(0,0){\line(0,1){4}}
 \put(4,4){\circle*{.3}}
 \put(0,4){\circle*{.3}}

 \put(-.5,- 0.8){$\delta _0$}
 \put(4.3,-.8){$\delta _\lambda $}
 \put(-.5,4.5){$\delta _\mu$}
 \put(4.3,4.5){$\delta _1$}

 \put(0,0){\line(1,1){4}}
 \end{picture}
              }
\put(7,0){ (iii)}
\put(22,0){ (iv)}

\end{picture}
\caption{(i) The paths $l_0$, $l_1$, $l_\lambda$, $l_\mu$ in the case when
$\lambda \not \in \R$ and $h_\lambda \neq h_\mu $;
(ii) Dynkin diagram in the case when $\lambda \not \in \R$ lies inside the
left branch of the hyperbola $\Gamma$;
(iii) Dynkin diagram in the case when $\lambda \not \in \R$ lies inside the
right branch of the hyperbola $\Gamma$;
(iv) Dynkin diagram in the case when $\lambda \not \in \R$ lies between the
two branches of the hyperbola $\Gamma$.}
\label{dynkin1}
\end{figure}

We finish this section by introducing a property of the continuous
families of ovals which is dependent on the structure of the related
Dynkin diagram. Suppose that the real polynomial $H= y^2+ P(x)$ has
distinct critical values. Let $\delta (h) \subset \{H=h\}$ be a
continuous family of ovals defined on a maximal open interval
$\Sigma =(h_c,h_s)$, where for $h=h_c$ the oval degenerates to a point
$\delta (h_c)$ which is a center and for $h=h_s$ the oval becomes a
homoclinic loop of the Hamiltonian system $dH=0$. The family
$\{ \delta (h)\}$ represents a continuous family of cycles vanishing
at the center $\delta (h_c)$.
To formulate our main result we shall need the following
\begin{definition}
\label{excd}
{\em It is said that $\{ \delta (h)\}$ is {\em an exceptional family of ovals},
provided that for every polynomial one-form $\omega $ the Abelian integral
$$
I(h)= \int_{\delta (h)} \omega , \quad h\in \Sigma
$$
allows an analytic continuation in the sector $S_\varepsilon (h_s)= \{h\in \C:
Arg (h-h_s) \in (-\varepsilon, 2\pi + \varepsilon )\}$ for some
strictly  positive
$\varepsilon$. The corresponding family of vanishing cycles represented by
$\delta (h)$ is called {\em an exceptional family of vanishing cycles.}}
\end{definition}
The above definition has in fact a geometric nature.
Indeed, $I(h)$ has an analytic continuation in
the sector  $S_\varepsilon (h_s)$ if and only if $\delta (h)$ has an
appropriate intersection number with each of the other families of
vanishing cycles. For instance, if $\lambda$, $\mu$ are real and
$h_1 < h_\mu < h_\lambda < h_0$, then the family of cycles $\delta_1(h)$ is
exceptional, because it has zero intersection number with $\delta _0(h)$.
On the contrary, the family $\delta _\mu (h)$ is not exceptional, because it
has a non-zero intersection number with  $\delta _0(h)$.
The remaining possible distributions of the critical
values, as well as the case $\lambda, \mu \not \in \R$, are studied in the
same way. We summarize all this in the following
\begin{proposition}
\label{exceptional}
Suppose that the Hamiltonian $H(x,y)$ is taken in a normal form
$(\ref{normal})$ and has four distinct critical values. Then:

\vspace{1ex}
\noindent
$(i)$
In the real case $(0 <\mu< \lambda <1)$, the continuous family of ovals
${\cal O}_1$ surrounding the  center
$(1,0)$ is exceptional, and the continuous family of ovals ${\cal O}_\mu$
surrounding the center $(\mu,0)$ is not exceptional.

\vspace{1ex}
\noindent
$(ii)$
In the complex case $(\lambda = \bar{\mu }\not \in \R)$,
the continuous family of  ovals ${\cal O}_1$
is exceptional if and only if $\lambda$ lies inside the left
branch of the hyperbola $\Gamma$ $($see
$(\ref{hyperbola})$ and {\em Fig. \ref{imconj})}.
\end{proposition}

\section{The period integral in a complex domain}
Assume that the level curve $\{H=h\}$ contains an oval $\delta(h)$ and
consider the integral (oriented along with the vector field (\ref{system}))
$$I_0(h)=\oint_{\delta(h)}\frac{dx}{y}$$
which is the derivative of the area inside $\delta(h)$ and hence
determines the period of the periodic orbit lying on $\delta(h)$.
Denote by $\alpha(h)$ and by $\beta(h)$ the minimal and the maximal
values of $x$ for $(x,y)\in\delta(h)$. These are among the real solutions
of the equation $P(x)=h$. Then one can express $I_0(h)$ in the
form
$$I_0(h)=\sqrt{2}\int_{\alpha(h)}^{\beta(h)}\frac{dx}{\sqrt{h-P(x)}}.$$
Consider first the case $0<\mu<\lambda<1$ and assume that $\delta(h)$ belongs
to either of the continuous families of
ovals ${\cal O}_\mu$, ${\cal O}_1$ surrounding centers $(\mu,0)$ and
$(1,0)$, respectively. The corresponding families of cycles vanish
respectively at $h=h_\mu$ and $h=h_1$. Denote this value in both cases
by $h_c$ and assume that the real ovals are defined in  $\Sigma=(h_c, h_s)$.
Finally, denote by $x_c$ the $x$-coordinate of the corresponding center.
We keep the same notation in the complex case and in the
cases when either $\mu=0$ or $\mu=\lambda$  or $\lambda=1$.
Then $x_c$ is the $x$-coordinate of the unique nondegenerate center
(as long as $(\lambda,\mu)\not\in\{(1,0), (1,1)\}$).

The integral $I_0(h)$ has an analytic continuation for $h<h_s$. We wish to
prove that $I_0(h)\neq 0$ in  $(-\infty, h_c)$. The precise statement
is as follows.

\begin{proposition}
\label{nozero}
Assume that the parameters $\lambda,\mu$ in $(\ref{normal})$ are either
real or lie inside the hyperbola $\Gamma$. Then the integral $I_0(h)$
related to any period annulus around a nondegenerate center of
$(\ref{system})$ has an analytic continuation in $(-\infty,h_s)$
and takes positive values there.
\end{proposition}
{\bf Proof.}
The claim that $I_0(h)$ has an analytic continuation in the interval
$(-\infty,h_s)$ is obvious in the complex case ($\lambda=
\bar{\mu}\not\in\R$). In the real case this follows from the Dynkin
diagram shown on Fig. \ref{realdynk}. For instance, if
$h_1 < h_\mu < h_\lambda < h_0$ as on Fig. \ref{realdynk}, and if
$h_s=h_\lambda$, $h_c= h_1$, then the claim is obvious. If $h_s=h_\lambda$,
$h_c= h_\mu$, then the claim is true again because the intersection number
of $\delta_\mu$ and $\delta_1$ is equal to zero. The remaining four
configurations of critical values and period annuli are examined in the
same way (with the help of Fig. \ref{realdynk}). Note, however, that the
integral $I_0(h)$ corresponding to the period annulus ${\cal O}_e$
(in this case $\Sigma =(h_\lambda, 0)$) has no analytic continuation on
the real axis outside $\Sigma$. This is so because the singularities in
both $h=h_0$ and $h=h_\lambda$ are unavoidable. The same conclusion holds
for some of the degenerate cases as well.

Next we shall prove the positivity of the period integral $I_0(h)$
on $(-\infty,h_s)$. We first choose a proper formula to present $I_0(h)$
when $h<h_c$. Take $z=x+iy\in\C$, $h<h_s$ and consider the polynomial
$h-P(z)$. By Taylor's formula,
$$\begin{array}{rl}
h-P(z)=&\!\!\!h-P(x)+\frac12y^2P''(x)-\frac{1}{24}y^4P''''(x)
-iy[P'(x)-\frac16y^2P'''(x)+\frac15y^4]\\[2mm]
\equiv&\!\!\!h-Q(x,y)-iyR(x,y)\end{array}$$
Given $h<h_s$, consider the equation  $h-P(z)=0$ and denote by
$\zeta_k(h)$, $k=1,2$, the two branches of the algebraic function defined by
$h-P(\zeta(h))\equiv 0$ which satisfy $\zeta(h_c)=(x_c,0)$.
These functions are unique since $P'(z)\neq 0$ for $z\neq 0,\mu,\lambda,1$.
For $h_c<h<h_s$, one has $\zeta_1(h)=(\alpha(h),0)=\alpha(h)$ and
$\zeta_2(h)=(\beta(h),0)=\beta(h)$ where $\alpha, \beta$ are as above.
For $h<h_c$, one has $\zeta_1(h)=\bar\zeta_2(h)\in\C$.
For a similarity in notation, we put $\zeta_1(h)=\alpha(h)$ and
$\zeta_2(h)=\beta(h)=\bar\alpha(h)$ where $Im\,\alpha<0$.
Denote by $C_h$ the curve in the complex plane connecting the
points $\alpha(h)$, $\beta(h)$ along $R(x,y)=0$, if $h<h_c$,
and along the real line $y=0$, if $h_c<h<h_s$. Then the integral
$I_0(h)$ is expressed as
\begin{equation}
\label{cmplx}
I_0(h)=\sqrt{2}\int_{C_h}\frac{dz}{\sqrt{h-Q(x,y)}},\quad h<h_s.
\end{equation}
Let $h<h_c$. We will establish below that $y=Im\,z$ can be used as a local
coordinate on $C_h$. Thus, $dz=[x'(y)+i]dy$ where $x(-y)=x(y)=x$. Denote for
short $y_h=Im\,\beta(h)$ and write $\int_{C_h}=\int_{-y_h}^0+\int_0^{y_h}$.
Replacing $y$ with $-y$ in the first integral, we obtain
\begin{equation}
\label{y}
I_0(h)=2\sqrt{2}\int_0^{y_h}\frac{dy}{\sqrt{Q(x(y),y)-h}},\quad h<h_c.
\end{equation}

Below we study the curve ${\cal R}:\;R(x,y)=0$ and establish that on this
curve the real-valued function $Q(x,y)-h$ is positive for all $h<h_c$.

Recall that $\cal R$ has an equation
$\frac15y^4-\frac16y^2P'''(x)+P'(x)=0$. We will study how this algebraic
curve changes when varying the parameters $\lambda,\mu$. It turns out that
a local analysis can be used to determine the global behavior of $\cal R$.
To begin with, we observe that the Poincar\'e index of any sufficiently big
circle $x^2+y^2=r^2$ subject to the vector field $dR(x,y)=0$ is $-3$.
The critical points of $R$ are determined from $R_x=R_y=0$, namely

\vspace{1ex}
(i) $y=0$, $P''(x)=0$;

\vspace{1ex}
(ii) $y^2=\frac{5}{12}P'''(x)$, $P''(x)-\frac{5}{72}P'''(x)P''''(x)=0$.

\vspace{1ex}
\noindent
It is easy to verify that the Hessian $R_{xx}R_{yy}-R_{xy}^2$ is negative
at any (real) nondegenerate critical point $(\xi,\eta)$. Indeed, it equals
$-\frac13(P'''(\xi))^2$ and
$-\frac{5}{108}P'''(\xi)(P''''(\xi))^2-\frac49(P'''(\xi))^2$, respectively
in cases (i) and (ii). Hence, $R$ has in general three nondegenerate saddles
and when two of the saddles collide, the resulting degenerate
critical point has only hyperbolic sectors around. In particular,
(ii) has no solutions for $0<\mu<\lambda<1$, since all the
three critical points are then given by (i). Also, in this case the curve
$\cal R$ is free of critical points because the value of $R$ at the saddle
$(\xi,0)$ is $P'(\xi)\neq 0$. In general, one can prove that $\cal R$
has no critical points outside the real line. Indeed, taking a critical point
$z=\xi+i\eta$ with $\eta\neq 0$, one has
$P'(z)=R(\xi,\eta)+\eta(R_y(\xi,\eta)+iR_x(\xi,\eta))=R(\xi,\eta)\neq 0$
unless $z=\lambda$ or $z=\mu$ (which may happen in the complex case).
However, by (ii), $(Re\,\lambda,\pm Im\,\lambda)$ is a critical point
if and only if $\lambda$ lies on the hyperbola $\Gamma$,
which is not our case. We have proved that the only critical points
$\cal R$  may have are $(0,0)$, $(\lambda,0)$ and $(1,0)$
provided that $\mu=0$, $\lambda=\mu$, $\lambda=1$, respectively.

As a consequence, the branches of the level curve $\cal R$ do not intersect
outside the real line and all they do escape to infinity (in the directions
of $x$ and $y$ altogether, as shown by the asymptotics at infinity). This
is because any closed compact curve necessarily surrounds a critical point
having an index $+1$ and $R(x,y)=const$ has no such points.

\begin{figure}[htbp]
\begin{center}
  \psfig{file=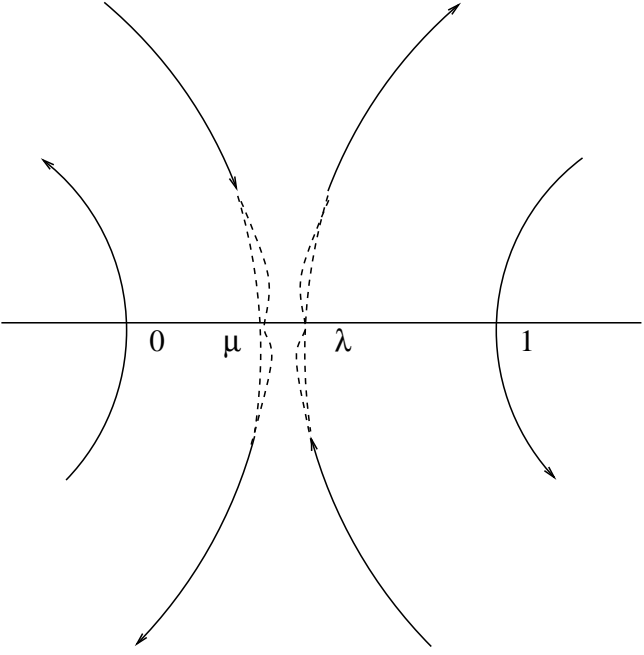}
\end{center}
\caption{The curve $R(x,y)=0$ for the  case $0<\mu<\lambda<1$.}
\label{curve}
\end{figure}

Assume first that $0<\mu<\lambda<1$. Then the four branches of $\cal R$ have
no common points at all. This implies that the branches through
$(0,0)$ and $(\mu,0)$ go to $-\infty$ in the $x$-direction while the
branches through $(\lambda,0)$ and $(1,0)$ go to $+\infty$, see Fig.
\ref{curve}. When running any of the branches, the $y$ variable changes
in a monotone way. (Otherwise, there would exist a horizontal line $y=const$
having at least 6 intersections with $\cal R$ which is impossible.)
Finally, denote by $\xi$ any of $0,1,\lambda,\mu$. The local equation
of $\cal R$ near the point $(\xi,0)$ is
$y^2=[6P''(\xi)/P'''(\xi)](x-\xi)+O((x-\xi)^2)$ where the coefficient
is negative for $\xi=0$, positive for $\xi=1$ and changes sign if
$\xi=\lambda,\mu$. This latter is because $P'''(\lambda)$ and
$P'''(\mu)$ change sign in $T$.
On the $(\lambda,\mu)$ curve where $P'''$ vanishes, the local equation
of $\cal R$
becomes $y^4\sim-5P''(\xi)(x-\xi)$ for $\xi=\lambda,\mu$.
Summing up, we obtain the form of the curve $\cal R$ as given in
Fig. \ref{curve}.
The arrows indicate the direction of the Hamiltonian vector field $dR=0$.

It is easy to see that $q(y)=Q(x(y),y)-h$ does not change sign when
$y\in(0,y_h)$. As $q(y_h)=0$, it suffices to show that $q'(y)\neq 0$.
We first note that $\dot{y}=-R_x<0$ on the branch through $(h_c,0)$
we consider. Then, by Cauchy-Riemann equations,
$$Q_x=(yR)_y,\quad Q_y=-(yR)_x,$$
and we obtain on $R=0$ that
$$q'=Q_y+x'Q_x=-\frac{y}{R_x}(R_x^2+R_y^2)<0.$$
This implies that $Q(x(y),y)-h>0$ for $y\in(0,y_h)$
and hence by (\ref{y}), $I_0(h)>0$ for $h<h_s$.

The same analysis is applicable when $\mu=0$ or $\mu=\lambda$ or $\lambda=1$,
the unique difference is that on Fig. \ref{curve} the corresponding
branches touch each other at $y=0$.

Let us turn now to the complex case. Apart of the real case, the branches
through $(0,0)$ and $(1,0)$ (let us denote them by $R_0$ and $R_1$) may
bifurcate significantly. This is because they are not necessarily monotone
with respect to $y$ and may escape to either $-\infty$ or $+\infty$, in the
$x$-direction. The complex plane is divided by the hyperbola $\Gamma$
into three parts. In the outer domain
$(Im\,\lambda)^2>5Re\,\lambda(Re\,\lambda-1)$, $R_0$ and $R_1$ behave as
in Fig. \ref{curve},
but are not necessarily monotone in $y$. Inside the left branch of the
hyperbola, $R_0$ and $R_1$ look like the branches through $\lambda$ and $1$
on Fig. \ref{curve}, respectively. And inside the right branch, $R_0$ and
$R_1$ behave as the branches through $0$ and $\mu$ on Fig. \ref{curve},
respectively. All this is easily verified by a simple deformation argument,
with a starting point on the appropriate side of triangle $T$.

To avoid complications, we regard only the parameters $\lambda$,
$\mu$ in the complex plane which belong to the interior of $\Gamma$
(which is the case to be considered in our main theorem). Taking
$\lambda=a+ib$, then $b^2<5a(a-1)$.
We are going to prove that $R_1$ is monotone with respect to $y$
in this domain. If so, one obtains just as in the real case that
$q(y)=Q(x(y),y)-h>0$ and the proof is complete.

Consider first the interior of the left branch of $\Gamma$. Then $R_1$ lies
in the half-plane $x\geq 1$. To establish monotonicity, we rewrite $R(x,y)$
as a polynomial in $t=x-1$. One obtains
\begin{equation}
\label{descartes}
R=t^4+(3-2a)t^3+(3-4a+a^2+b^2-2y^2)t^2+(1-2a+a^2+b^2+(2a-3)y^2)t+R(1,y).
\end{equation}
Assume that $a\leq\frac12$. As the first two coefficients of the polynomial
are positive, then according to Descartes rule of signs, (\ref{descartes})
can have three positive roots only if the coefficient at $t^2$ is negative
and the coefficient at $t$ is positive. However, these two conditions
contradict each other. Therefore each line $y=const$ can intersect $R_1$
in at most two points (counting multiplicity) and hence in a single point
because the number of intersections should be an odd number.

Now take $\lambda$ in the interior of the right branch of $\Gamma$ and denote
by $(\xi,\eta)$ the right-most point of $R_1$ in the half-plane $y\geq 0$.
Thus, $\xi\geq 1$ and $\eta\geq 0$. As above, one obtains with
$t=\xi-x \geq 0$ that
\begin{equation}
\label{descart1}
{\textstyle
R=t^4-\frac16P''''(\xi)t^3+(\frac12P'''(\xi)-2y^2)t^2-(P''(\xi)
-\frac16P''''(\xi)y^2)t+\frac15(y^2-\frac{5}{12}P'''(\xi))^2.}
\end{equation}
As $\eta=0$ if $P'''(\xi)\leq 0$ and $\eta^2=\frac{5}{12}P'''(\xi)$
otherwise, the coefficient at $t^2$ is negative for $y>\eta$.
Then by Descartes rule,
(\ref{descart1}) has only two positive roots and therefore the part of
$R_1$ above the line $y=\eta$ is monotone. It remains to consider the
part below it (if $\eta>0)$. In this domain (\ref{descart1}) can
have two complex roots and Descartes rule does not yield the required result.
We can proceed as follows. Assume that the part of $R_1$ in the upper
half-plane lying below $y=\eta$
is not monotone. Denote by $(\xi_1,\eta_1)$ the first local maximum on $R_1$
(starting from $(1,0)$), thus $1<\xi_1<\xi$ and $0<\eta_1<\eta$. Consider
also the branch $R_+$ of $\cal R$ which is placed in the upper half-plane
to the right of the line $x=\xi$ (and above the line $y=\eta_1$).
An easy calculation yields the equations
$y=\pm\sqrt{5\pm2\sqrt5}(x-\frac{1+2a}{4})$
for the tangential lines of $\cal R$ at infinity. As their common point
is different from $(\xi_1,\eta_1)$, there would always exist a tangential
line $l$ to $R_+$ at a finite point, which goes through $(\xi_1,\eta_1)$.
However, then $l$ would intersect $R_+$ in at least two points and $R_1$
in at least three points, which is impossible because $\cal R$ is a quartic
curve. The monotonicity of $R_1$ is established. $\Box$

\vspace{2ex}
\noindent
\begin{proposition}
\label{nozero1}
Assume that $\delta(h)$ is an exceptional family of ovals according to
Definition $\ref{excd}$. Then the period integral
$I_0(h)= \int_{\delta(h)} dx/y$
has an analytic continuation in the complex domain ${\cal D}= \C \setminus
[h_s,\infty)$ and does not vanish there. \end{proposition}
{\bf Proof.}
The proof that $I_0(h)$ has an analytic continuation in $\cal D$
follows immediately from the definition of an exceptional family of ovals.
To count the zeros of the analytic function $I_0(h)$ in $\cal D$,
we apply the argument principle to the complex domain
${\cal D}_r$  obtained from $\cal D$ by removing the discs
$\{h\in \C: |h-h_{s}|< r \}$, $\{h\in \C: |h| > 1/r \}$, where $ r>0$
is a small enough constant.  For this purpose we
evaluate the increase of the argument of $I_{0}(h)$ when $h$ runs
the boundary of ${\cal D}_r$ in a positive direction.

We consider first the real case ($0< \mu < \lambda < 1$).
The Picard-Lefschetz formula implies that in a suitable complex
neighborhood of $h_s$ it holds
$$I_{0}(h) = \varphi(h) \log(h-h_{s}) + \psi(h), \quad
\varphi(h)= \frac{1}{2\pi} \int_{\delta_{s}(h)} \frac{dx}{y}
$$
where $\varphi,\psi$ are locally holomorphic and, as it is easily
checked, $\varphi(h_s)\neq 0$.

The behavior of $I_{0}(h)$ for $|h|$ sufficiently big follows from
the Picard-Lefschetz formula too, but the argument is more delicate.
Let $l$ be a simple closed path which makes one turn in a positive
direction around all the critical values of $H$. It induces a
monodromy map
$$
l_* : H_1(\Gamma_{\tilde{h}},\Z) \rightarrow H_1(\Gamma_{\tilde{h}},\Z),
\quad \Gamma_{\tilde{h}}= \{(x,y)\in \C:\, H(x,y)=\tilde{h} \} \; .
$$
It is shown for instance in \cite{gav97} that $l_{*}$ coincides
with the operator of classical monodromy of the singularity
$\frac{1}{2}y^2 + \frac{1}{5}x^5 $ (which is the highest
weighted homogeneous part of $H$). The characteristic polynomial $P(z)$
of $l_*$ is now easily computed to be equal to
$P(z)= (z^5+1)/(z+1)$ (see \cite{bri}). This shows that in a neighborhood
of $h=\infty$ on the projective sphere the integral
$I_0(h)$ is a meromorphic function of $h^{1/10}$. The substitution
$$
x \rightarrow x h^{1/5},\quad y \rightarrow y h^{1/2},\quad h \rightarrow
\infty
$$
defines an isotopy of the regular fibers of $H$ and the regular fibers
of its weighted homogeneous part
$\frac{1}{2}y^2 + \frac{1}{5}x^5 $ \cite{gav97}.
From this fact we deduce that
$$I_{0}(h)=c_{\infty}h^{1/5-1/2} (1+O(h^{-1/10})).$$
For a further use we note that
$c_{\infty} \neq 0$. Indeed, the one-form $dx/y$ is not cohomologous to zero
on the regular fibers of the polynomial $\frac{1}{2}y^2 + \frac{1}{5}x^5 $
(this follows from \cite[Theorem 1.1]{gav98a}).
On the other hand, the form of the characteristic polynomial of $l_{*}$ shows
that it is irreducible over the field of rational numbers. Therefore the span
of $l_*^k \delta$, $k=0, 1, 2, 3$ generates the whole
homology group of the fiber. This implies that if $c_\infty=0$, then the
restriction of  $dx/y$ on the regular fibers of
$\frac{1}{2}y^2 + \frac{1}{5}x^5$  defines the zero
cohomology class which is a contradiction.

The above consideration shows that the decrease of the argument of $I_{0}(h)$
along $\{h\in \C: |h| = 1/r \}$ is close to $3 \pi/5$ and the increase of the
argument of $I_{0}(h)$  along $\{h\in \C: |h-h_{s}| = r \}$ is close to zero.

We claim further that the imaginary part of $I_{0}(h)$ along
$(r+h_s,1/r)$ does not vanish. Namely, denote the two determinations of
$I_{0}(h)$ on $(r+h_s,1/r)$ by $I_{0}^{\pm}(h)$. Similarly, denote by
$\delta^{\pm}(h)$ the two determinations of the continuous family
$\delta(h)$ on $(r+h_s,1/r)$.
Since $I_{0}$ is a real analytic function  along $(-\infty, h_{s})$, we have
$$I_{0}^{+}(h)= \overline{I_{0}^{-}(h)},\quad
\delta^{+}(h) =\overline{\delta^{-}(h)}, \quad h>h_{s}
$$
and hence $2\,\sqrt{-1} \, Im\, I_{0}^{+}(h) =
I_{0}^{+}(h) - I_{0}^{-}(h)$. On the other hand $\delta(h)$ is an exceptional
family of ovals. Therefore the Picard-Lefschetz formula implies that for
$h> h_s$,
$$
\delta^{+}(h) - \delta^{-}(h) = \gamma(h),\qquad
I_{0}^{+}(h) - I_{0}^{-}(h) = \int_{\gamma(h)} \frac{dx}{y}
$$
where $\gamma(h)$ is a continuous family of cycles vanishing at
$h=h_{s}$, with intersection number $\langle\gamma,\delta^-\rangle=1$.

Finally, we note that after substituting $y$ by $\sqrt{-1} y$,
the family of cycles $\gamma(h)$ corresponds to a family
of ovals surrounding a center. Moreover, in the real case the full
bifurcation diagram (consisting of the boundary of $T$ and the parts of
the curves $h_\lambda=h_0$, $h_\mu=h_1$, $h_0=h_1$ inside) images onto
itself. For later use we also mention that in the complex case, the
interior of the hyperbola $\Gamma$ images onto itself. All this holds
because the corresponding Hamiltonians when taken in a normal form
(\ref{normal}) are related by the formula
$$H(x, \sqrt{-1}y,\lambda,\mu)=h_1-H(1-x,y,1-\mu,1-\lambda).$$
Hence, it follows from Proposition \ref{nozero}  that
$\int_{\gamma(h)} dx/y$ does not vanish along $(h_{s},\infty)$. Therefore
if  $h\in (r+h_s,1/r)$, then $Im\, I_{0}^{\pm}(h) \neq 0$ and hence
the increase of the argument of $Im\, I_{0}^{\pm}(h)$
along $(r+h_{s},1/r)$ is at most $\pi$.

Summing up the above information we conclude that the increase of the
argument of $I_{0}(h)$ along the boundary of ${\cal D}_r$ is strictly
less than $2\pi$. By the argument principle, $I_{0}(h)$ does not vanish
in ${\cal D}_r$, and hence in $\cal D$.

The complex case when $\lambda = \bar{\mu }\not \in \R$ is studied in the
same way. $\Box$

\section{The Chebyshev property}
Suppose that the real polynomial $H= y^2+ P(x)$ has distinct critical values.
Let $\delta (h) \subset \{H=h\}$ be a continuous family of ovals defined on
a maximal open interval
$\Sigma =(h_c,h_s)$, where for $h=h_c$ the oval degenerates to a point
$\delta (h_c)$ which is a
center and for $h=h_s$ the oval becomes a homoclinic loop of the
Hamiltonian system $dH=0$. Denote ${\cal D}=\C\setminus[h_s,\infty)$.
The following theorem is our main result.

\begin{theorem}
\label{main}
Let $\delta (h)$, $h\in \Sigma $ be an exceptional family of ovals
of $(\ref{normal})$. Then the real vector space
${\cal A}=\{\alpha_0 I_0 + \alpha_1 I_1: \alpha_0, \alpha_1 \in \R \}$
of Abelian integrals
$$I_0(h) =  \int_{\delta (h)} \frac{ dx}{y}, \qquad I_1(h) =
\int_{\delta (h)} \frac{x dx}{y}, \qquad h\in \Sigma $$
is Chebyshev in $\cal D$ and hence in $\Sigma$.
\end{theorem}
{\bf Proof.}
Let $\delta(h)$ where $h \in \Sigma = (h_{c},h_{s})$ be an exceptional
family of ovals. According to Proposition \ref{nozero1}, it suffices to
show that the analytic function $F(h)= \alpha_0 + I_{1}(h)/I_{0}(h)$ has
at most one zero in the domain $\cal D$. To prove this claim we apply the
argument principle to $F(h)$ in the complex domain ${\cal D}_r$ obtained
from $\cal D$ by removing the discs
$\{h\in \C: |h-h_{s}|< r \}$, $\{h\in \C: |h| > 1/r \}$,
where $ r>0$ is a small enough constant.

 As in the proof of Proposition \ref{nozero1}, we may check that:
\begin{itemize}
     \item  if $h\sim \infty$, then $I_{1}(h)$ is meromorphic
      with respect to $h^{-1/10}$ and
      $$I_{1}(h) = c_{\infty} h^{2/5-1/2}(1+O(h^{-1/10}))
      \;\;\mbox{\rm where}\;\; c_{\infty} \neq 0;$$

      \item if $h \sim h_{s}$, then
      $$I_{1}(h) = \varphi_1(h) \log(h-h_{s}) + \psi_1(h),
 \quad \varphi_1(h)= \frac{1}{2\pi} \int_{\delta_{s}(h)} \frac{x\,dx}{y},$$
  where $\varphi_1$ and $\psi_1$ are locally holomorphic. Moreover,
  it is easily verified that
  $\lim\limits_{h \rightarrow h_{s}} I_{1}(h)/I_{0}(h)=
  \varphi_1(h_s)/\varphi(h_s) = x_{s} $
  where $(x_{s},0)$ is the saddle point corresponding to the
  critical value $h_{s}$.
\end{itemize}

\noindent
Therefore the increase of the argument of $F(h)$ along
$\{h\in \C: |h| = 1/r \}$
is close to $2 \pi/5$ and the increase of the argument of $F(h)$
along $\{h\in \C: |h| = r \}$ is close to zero or strictly negative.

We claim further that
\begin{itemize}
        \item  the imaginary part of $F(h)$ along $(r+h_{s},1/r)$
        does not vanish.
\end{itemize}

\noindent
If this claim were proved then summing up the above information we would get
that the argument of $F(h)$ increases along the boundary of ${\cal D}_r$
by at most $2\pi + 2 \pi/5$. Thus the
function $F(h)$ can have at most one zero in  ${\cal D}_r$  and hence
in $\C \backslash [h_{s}, \infty)$.

It remains to prove the above claim. Recall that, since $\delta (h)$ is an
exceptional family of ovals, then the two determinations of $F(h)$
over $(h_s,\infty)$, namely $F^\pm(h)$, are analytic
along $(h_s,\infty)$. The function $F(h)$ is real analytic along
$(-\infty,h_s)$ and hence
$$
2\,\sqrt{-1}\,Im\, F(h) = F^+(h) - F^-(h) \;\; \forall h\in (h_s,\infty).
$$
Denote by $\delta ^\pm(h)$ the two determinations of
$\delta (h)$, $h\in {\cal D}$ on $(h_s,\infty)$, and let
$$\Gamma_{h}~=~\{(x,y)\in \C^2: H(x,y) =h\}.$$
It was shown in the proof of Proposition~\ref{nozero1} that
$$
\int_{\delta^{\pm} (h)} \frac{ dx}{y} \neq 0\quad\forall h \in (h_s,\infty).
$$
Then we obtain for each $h\in (h_s,\infty)$
$$
   2\,\sqrt{-1}\,Im\, F(h)  =  \frac{\int_{\delta^+(h)} \frac{x dx}{y}}
   {\int_{\delta^{+} (h)} \frac{dx}{y}}
                  - \frac{\int_{\delta^{-} (h)} \frac{x dx}{y}}
   {\int_{\delta^{-} (h)} \frac{ dx}{y}}
 =  \frac{ det \left( \begin{array}{cc}
                \int_{\delta^{+} (h)} \frac{x dx}{y} & \int_{\delta^{+}
(h)} \frac{ dx}{y}  \\[1mm]
                \int_{\delta^{-} (h)} \frac{x dx}{y} & \int_{\delta^{-}
(h)} \frac{ dx}{y} \end{array} \right) }{
        |\int_{\delta^{\pm} (h)} \frac{ dx}{y}|^{2}}.
$$
Therefore we have to prove that the analytic function
$$
\Delta (h) =  det \left( \begin{array}{cc}
                \int_{\delta^{+} (h)} \frac{x dx}{y} & \int_{\delta^{+}
(h)} \frac{ dx}{y}  \\[1mm]
                \int_{\delta^{-} (h)} \frac{x dx}{y} & \int_{\delta^{-}
(h)} \frac{ dx}{y} \end{array} \right),\quad  h \in (h_s,\infty)
$$
does not vanish.

Suppose first that $h=h_0 \in (h_s,\infty)$ is a critical value of $H$.
The affine curve $\Gamma_{h_0}$ is singular of arithmetic genus two. Let
$\bar{\Gamma} _{h_0}$ be the corresponding completed and normalized Riemann
surface. As the geometric genus of $\bar{\Gamma} _{h_0}$ is one and the intersection
number $\langle \delta^{+},\delta^- \rangle$ equals $+1$, then
$\delta^{+}(h_0)$ and $\delta^-(h_0)$ form a canonical homology basis of
$\bar{\Gamma} _{h_0}$. We have $P(x)=(x-x_0)^2(x-x_1)(x-x_2)(x-x_3)$
where $x_i\neq x_j$ for $i\neq j$. Therefore, $(x-x_0)dx/y$ induces a
holomorphic one-form on $\bar{\Gamma} _{h_0}$ and $dx/y$ induces a
meromorphic one-form of the third kind on $\bar{\Gamma} _{h_0}$ with simple
poles at $P^\pm$, where $P^\pm$ are the two distinct pre-images
of the singular point $(x_0,0)\in \Gamma _{h_0}$ under the canonical
projection $ \bar{\Gamma} _{h_0}\rightarrow \Gamma_{h_0}$. The reciprocity
law for differentials of the first and third kind \cite{gh} applied to
$(x-x_0)dx/y$, $dx/y$ on the elliptic curve $\bar{\Gamma} _{h_0}$ gives
\begin{eqnarray*}
\Delta (h) & = & 2\pi \sqrt{-1}\left[
Res|_{P^+}\left(\frac{dx}{y} \int_{P^0}^{P} \frac{(x-x_0)dx}{y}\right)
+ Res|_{P^-}\left(\frac{dx}{y} \int_{P^0}^{P} \frac{(x-x_0)dx}{y}\right)
\right]\\
&=& \frac{\pm 2 \pi \sqrt{-1}}{\sqrt{(x_0-x_1)(x_0-x_2)(x_0-x_3)}}
\int_{P^-}^{P^+}\frac{dx}{\sqrt{(x-x_1)(x-x_2)(x-x_3)}
}
\end{eqnarray*}
where the path of integration is on the Riemann surface
$ \bar{\Gamma} _{h_0}$ of the affine elliptic curve
$$\{y^2= (x-x_1)(x-x_2)(x-x_3) \}$$
and $P^\pm = (x_0,\pm \sqrt{(x_0-x_1)(x_0-x_2)(x_0-x_3)}) $. The function
$$
P \rightarrow \int^P \frac{dx}{\sqrt{(x-x_1)(x-x_2)(x-x_3)}}
$$
is an uniformizing variable on the elliptic curve $ \bar{\Gamma} _{h_0}$
which shows that
$$
\int_{P^-}^{P^+}\frac{dx}{\sqrt{(x-x_1)(x-x_2)(x-x_3)}} = 0
\;\;\Rightarrow\;\; P^+=P^- .
$$
As $P^+ \neq P^-$, then $\Delta (h_0) \neq 0$.

Suppose now that $h\in (h_s,\infty)$ is a regular value of $H$. The
intersection number $\langle \delta^{+},\delta^- \rangle$ equals $+1$ and
it is easy to check that we may choose a canonical homology basis
$\delta_{1}(h),\ldots,\delta_{4}(h)$ of the lattice
$H_{1}(\bar{\Gamma}_{h},\Z)$ such that $\delta_{1}(h) = \delta^{+}$,
$\delta_{3}(h) = \delta^{-}$. Further, let $\omega_{1}, \omega_{2}$ be a
normalized basis of holomorphic differentials on $\bar{\Gamma}_{h}$. The
corresponding period lattice $\Pi$ reads
$$
\Pi = \left(
\begin{array}{cccc}
        \int_{\delta_{1} (h)} \omega_{1} & \int_{\delta_{2} (h)} \omega_{1} &
        \int_{\delta_{3} (h)} \omega_{1} & \int_{\delta_{4} (h)} \omega_{1}\\
\int_{\delta_{1} (h)} \omega_{2} & \int_{\delta_{2} (h)} \omega_{2} &
        \int_{\delta_{3} (h)} \omega_{2} & \int_{\delta_{4} (h)} \omega_{2}
\end{array}
\right) =
\left(
\begin{array}{cccc}
        1 & 0 &
        a & b  \\
0 & 1 &
        c & d
\end{array}
\right)
$$
where
$$
\left(
\begin{array}{c}
        \omega _1 \\
        \omega _2
\end{array}
\right) =
\left(
\begin{array}{cc}
   \int_{\delta_{1} (h)} \frac{xdx}{y} & \int_{\delta_{2} (h)}\frac{xdx}{y}\\
   \int_{\delta_{1} (h)} \frac{dx}{y} & \int_{\delta_{2} (h)}\frac{dx}{y}
\end{array}
\right)^{-1}
\left(
\begin{array}{c}
        \frac{xdx}{y} \\
      \frac{ dx}{y}
\end{array}
\right) .
$$
We get that
$$
\Delta (h) = \det \left(
\begin{array}{cc}
    \int_{\delta_{1} (h)} \frac{xdx}{y} & \int_{\delta_{2} (h)}\frac{xdx}{y}\\
    \int_{\delta_{1} (h)} \frac{dx}{y} & \int_{\delta_{2} (h)}\frac{dx}{y}
\end{array}
\right) \det
\left(
\begin{array}{cc}
        1 & a \\
     0 & c
\end{array}
\right)
$$
and hence
$$
\Delta (h)=0 \;\;\Leftrightarrow \;\; c=0.
$$
Let $\Lambda $ be the lattice generated by the columns of $\Pi $.
The Riemann bilinear relations \cite{gh} on $\Pi $ imply
$$
b=c.
$$
Therefore if $\Delta (h)=0$, then the Jacobian variety
$J(\bar{\Gamma} _{h}) = \C^2/\Lambda $ is
a direct product of the two elliptic curves
\begin{equation}
\label{curves}
\C / \{\Z \oplus a \Z  \},\quad \C / \{\Z \oplus d \Z  \}.
\end{equation}
It is well known that this is impossible (e.g. the Remark on page 49 in
\cite{arbarello}) which completes the proof of Theorem \ref{main}.

For convenience of the reader, we present a proof of the last claim.
Consider the Riemann theta function
$$
\theta(z)=  \sum_{m\in \Z^{2}} e^{\pi \sqrt{-1} \langle m,Bm\rangle + 2 \pi
\sqrt{-1}\langle m,z\rangle}
$$
and the corresponding theta divisor
$\Theta = \{z \in  J(\bar{\Gamma}_{h}): \;\theta (z)=0 \}$ where
$$
B= \pmatrix{ a & b  \cr c & d},\quad z= \pmatrix{z_{1}\cr z_{2}},\quad
 \langle m,z\rangle = m_{1}z_{1}+ m_{2}z_{2}.
$$
By the Riemann theorem, the divisor $\Theta$ is isomorphic to the curve
$\bar{\Gamma} _{h}$ and in particular is irreducible. If $\Delta(h)=0$,
then $b=c=0$ and $\theta(z)$ factorizes into a product of two elliptic
theta functions
$$
\theta(z)=
\sum_{m_{1}\in \Z} e^{\pi \sqrt{-1} a m_{1}^{2} + 2 \pi
\sqrt{-1} m_{1} z_{1}}
\sum_{m_{2}\in \Z} e^{\pi \sqrt{-1} d m_{2}^{2} + 2 \pi
\sqrt{-1} m_{2} z_{2}}  .
$$
Therefore the Riemann theta divisor $\Theta$ is a product of the two
elliptic curves (\ref{curves}) which is a contradiction.  $\Box$

\section{Concluding remarks}
We finish this paper with some general hypotheses concerning the zeros of
the function $\alpha_0I_0+\alpha_1I_1$. In the previous section we
demonstrated that this linear combination belongs to a Chebyshev space,
provided that the continuous family of ovals $\{\delta(h)\}$ is exceptional.
In the non-exceptional case a big part of our proof still holds,
provided that $\lambda$ and $\mu$ are real. This fact combined
with the exact local result obtained in Proposition \ref{cyclicity}
suggests the following more general statement.

\vspace{2ex}
\noindent
{\bf Conjecture 1.}
Assume that the polynomial $P(x)$ has distinct real critical values.
 Then the corresponding
vector space of Abelian integrals $\alpha_0 I_0+\alpha_1I_1$
is Chebyshev with accuracy at most one in ${\cal D}$.


\vspace{2ex}
\noindent
One may also ask whether the upper bound (two) for the number of
the zeros near the center always holds in the whole interval $\Sigma$ and
next, whether the bounds for the real zeros when $h\in\Sigma$ and for the
complex zeros when $h$ belongs to a certain complex extension of $\Sigma$
do coincide. It turns out that the answers of all these questions are
negative. In particular, we establish below that when $P(x)$ has
complex critical values, then there are three zeros in $\Sigma$
for appropriate $\lambda\in\C$.

One can try to predict the possible number of real zeros by studying
the behavior of the ratio $F(h)=I_1(h)/I_0(h)$ near the ends of
the interval $\Sigma$. For reader's convenience, we first recall some
formulas.

Let us take a continuous function $f(x)$ and consider the integral
$$I(h)=\int_{\delta(h)}\frac{f(x)dx}{y},\quad h\in\Sigma$$
where $\Sigma=(h_c,h_s)$ for the families ${\cal O}_1$, ${\cal O}_\mu$
and $\Sigma=(h_s,h_0)$ for the family ${\cal O}_e$.
Using Picard-Lefschetz formula or direct calculations of integrals,
one can obtain the following asymptotical expansions in a neighborhood
of the critical levels $h_s$ corresponding to a saddle or a cuspidal point
$(x_s,0)$:

\vspace{1ex}
\noindent
(i) For the saddle-loop and eight-loop cases,
$$I(h)=-c_sf(x_s)\log|h-h_s|+\varphi(h)$$
where $c_s$ is a specific positive constant independent on $f$ and
$\varphi(h)$ is a continuous function in a neighborhood of $h_s$.
We are not going to write up the exact expression of $c_s$ since do
not need it here.

\vspace{1ex}
\noindent
(ii) For the heteroclinic loop cases,
$$I(h)=-(c_0f(0)+c_s f(x_s))\log|h-h_s|+\varphi(h)$$
where $c_0$, $c_s$ and $\varphi(h)$ are as above.

\vspace{1ex}
\noindent
(iii) For the cuspidal loop cases,
$$I(h)=c_sf(x_s)|h-h_s|^{-1/6}+\varphi(h),$$
where $c_s$ and $\varphi(h)$ are as above.

\vspace{1ex}
\noindent
(In fact, these formulas hold in a much more general context.)
Applying the above expansions to $F(h)$, we get $F(h_s)= x_s/(1+c_0/c_s)$
for the heteroclinic loop case and $F(h_s)=x_s$ for all the remaining cases.
Moreover, what concerns the saddle-loops, it is evident that
$F(h)-x_s\neq 0$ in $\Sigma$ and the sign of the expression depends on
the position of the saddle with respect to the family of ovals.
Consider now the case of ${\cal O}_e$ where the situation is more
complicated. For $h=h_\lambda$, we obtain from (\ref{normal}) that
$$\textstyle y^2=\frac25(x-\lambda)^2p(x),\quad p(x)=(x-x_1)(x-x_2)(x_3-x)$$
where $x_1<0<x_2\leq\lambda\leq x_3$ and $x_2$, $x_3$ are the left-most and
the right-most points of the eight-loop (or the cuspidal loop, as a
limit case). By (i)--(iii), the function $\varphi(h)=I_1(h)-\lambda I_0(h)$
is continuous in a neighborhood of $h_\lambda$, which yields immediately that
$$\varphi(h_\lambda)=
\sqrt{10}\left(\int_\lambda^{x_3}-\int_{x_2}^\lambda\right)\frac{dx}
{\sqrt{p(x)}}.$$
Clearly, $\varphi(h_\lambda)\to-\infty$ when $x_2\to x_1$,
$\varphi(h_\lambda)\to\pm \sqrt{10}\int_{x_2}^{x_3}dx/\sqrt{p(x)}$
when $\lambda\to x_2$ and $\lambda\to x_3$, respectively.
As $x_2=x_1$ is equivalent to $h_\lambda=0$,
$\lambda=x_2$ is equivalent to $\lambda=\mu$ and $\lambda=x_3$
is equivalent to $\lambda=1$, this means that the curve
$$\gamma_s=\{(\lambda,\mu)\in T:\;I_1(h_\lambda)-\lambda I_0(h_\lambda)=0\}$$
connects the points $(0,0)$ and $(1,1)$ and lies in the part of $T$ above
$\gamma$. Surprisingly, although expressed in terms of incomplete elliptic
integrals, it turns out that $\gamma_s$ is a part of an algebraic curve of
degree 6. This can be verified by transforming the integrals to a standard
form. From the explicit algebraic equation thus obtained it is easily seen
that $\gamma_s$ is a simple and connected curve (we omit the details).
The curve $\gamma_s$ is an element of the bifurcation diagram because for
$h$ close to $h_\lambda$, the function $F(h)-\lambda$ will change the sign
when crossing $\gamma_s$ in the parameter space.

Among other, the above analysis implies that for $\lambda=\mu$, the function
$\alpha_0I_0+\alpha_1I_1=I_0(\alpha_0+\alpha_1 F)$ can have three zeros
in $(h_1,h_\lambda)\cup(h_\lambda,h_0)$. Hence, the same fact remains true
for ${\cal O}_1$ and the respective $\Sigma=(h_1,h_0)$ in the complex case,
at least for parameters $\lambda=a+ib$ with $a\in(0,1)$ and $|b|$ small.
Based on the behavior of $F$ near the other endpoint which corresponds
to a center (see below), we conjecture that the domain in the parameter space
$\lambda\in\C$ where the equation $F(h)=const$ can have three real zeros,
is surrounded by a curve on which they join into a triple zero.

Looking at the bifurcation diagram again, let us denote by $\gamma_c$
the curve
$$\textstyle\gamma_c=\{(\lambda,\mu)\in T:\;
\lambda=\frac{3\mu^2-2\mu}{2\mu-1}\}
\quad\mbox{\rm (in the real case),}$$
$$\gamma_c=\{\lambda\in \C:\; |\lambda-2|=1\}
\quad\mbox{\rm (in the complex case).}$$
Further, let us denote by $\gamma^*$ the (hypothetical) curve
$$\gamma^*=\{\lambda\in \C:\;
\exists h^*\in\Sigma\;\; \mbox{\rm with}\;\; F'(h^*)=F''(h^*)=0\}.$$
We believe that $\gamma^*$ is a simple closed curve in the complex plane
going through $(0,0)$ and $(1,0)$. On $\gamma^*$, a bifurcation
of a triple real zero into a simple one occurs.

Next, in the real case, denote by $\Omega_\mu$ the part of $T$ placed
between the curves $\gamma_c$ and $\gamma$, by $\Omega_e$ the part of
$T$ placed between the curve $\gamma_s$ and the line $\lambda=\mu$
and finally, in the complex case, denote by $\Omega_1$
the unit disk inside $\gamma_c$ and by $\Omega_1^*$
the interior of $\gamma^*$ (see figures \ref{imconj} and  \ref{reconj}).

        \begin{figure}
\setlength{\unitlength}{.04\textwidth}

\begin{picture}(25,15)(5,0)
\put(2,8){
 \put(5,0){\vector(1,0){20}}
\put(12,-7){\vector(0,1){15}}
\qbezier[200](9,-6,)(14.96,0)(9,6)
\put(4,0){\qbezier[200](15,-6,)(9.04,0)(15,6)}
\put(12,0){
\qbezier[50](0,0)(0,6)(2,6)
\qbezier[50](0,0)(0,-6)(2,-6)
\qbezier[50](2,6)(4,6)(4,0)
\qbezier[50](2,-6)(4,-6)(4,0)
}
\put(16,0){
\qbezier[100](0,0)(0,2)(2,2)
\qbezier[100](0,0)(0,-2)(2,-2)
\qbezier[100](2,2)(4,2)(4,0)
\qbezier[100](2,-2)(4,-2)(4,0)}
\put(25.5,-.3){$Re \lambda$}
\put(13.5,.3){$\Omega_1^*$}
\put(17,.3){$\Omega_1$}
\put(11,-.8){$0$}
\put(16.5,-.8){$1$}
\put(20.5,-.8){$3$}
\put(15,-6){$\gamma^*$}
\put(19.5,-2){$\gamma_c$}
\put(19.7,6){$\Gamma$}
\put(8,6){$\Gamma$}
\put(11,8.5){$Im \lambda$} }
\end{picture}
\caption{Partition of the parameter space in the complex case 
$\lambda = \bar{\mu} \in \C$.}
\label{imconj}
\end{figure}
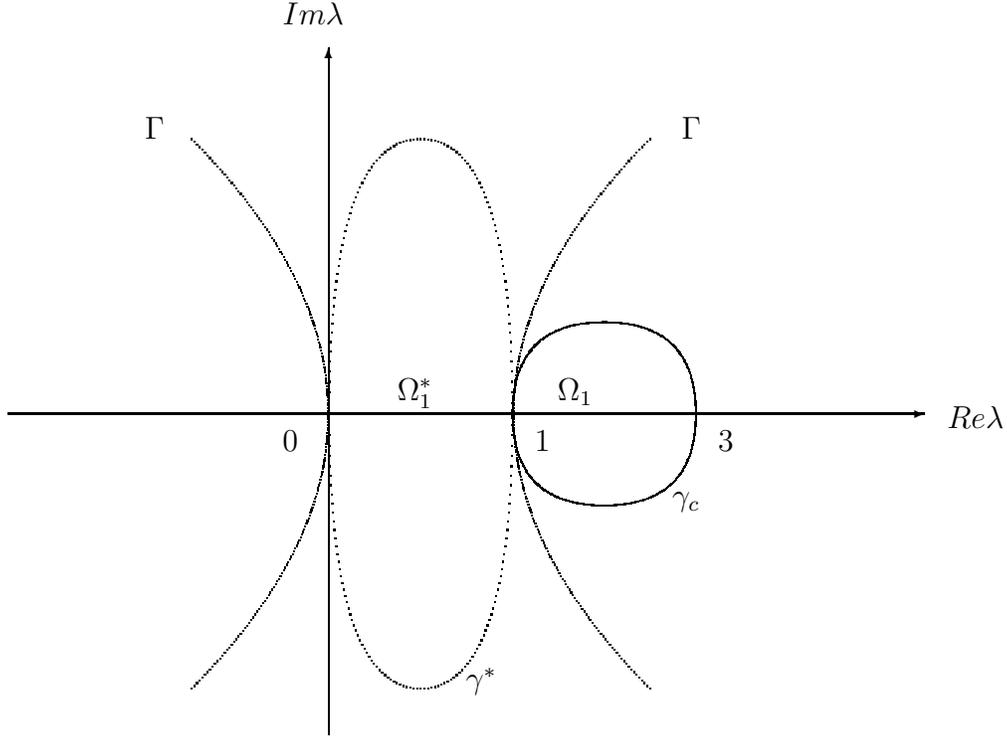
\begin{figure}
\setlength{\unitlength}{.04\textwidth}

\begin{picture}(25,25)(0,0)
\put(0,0){\vector(1,0){22}}
\put(0,0){\vector(0,1){20}}
\put(20,0){\line(0,1){20}}

\qbezier[100](0,0)(10,4.5)(20,6)
\qbezier[100](0,0)(9,6)(20,10)
\qbezier[100](0,0)(15,10)(20,20)
\put(0,0){\line(1,1){20}}
\put(0,-.6){$0$}
\put(20,-.6){$1$}
\put(22,.6){$\lambda$}
\put(0.6,19.6){$\mu$}
\put(10,3.2){$\gamma_c$}
\put(16,7){$\Omega_\mu$}
\put(14,8.1){$\gamma$}
\put(17,14.6){$\gamma_s$}
\put(10,9.1){$\Omega_e$}

\end{picture}
\caption{Partition of the parameter space in the real case $\lambda,\mu \in \R$.
}
\label{reconj}
\end{figure}
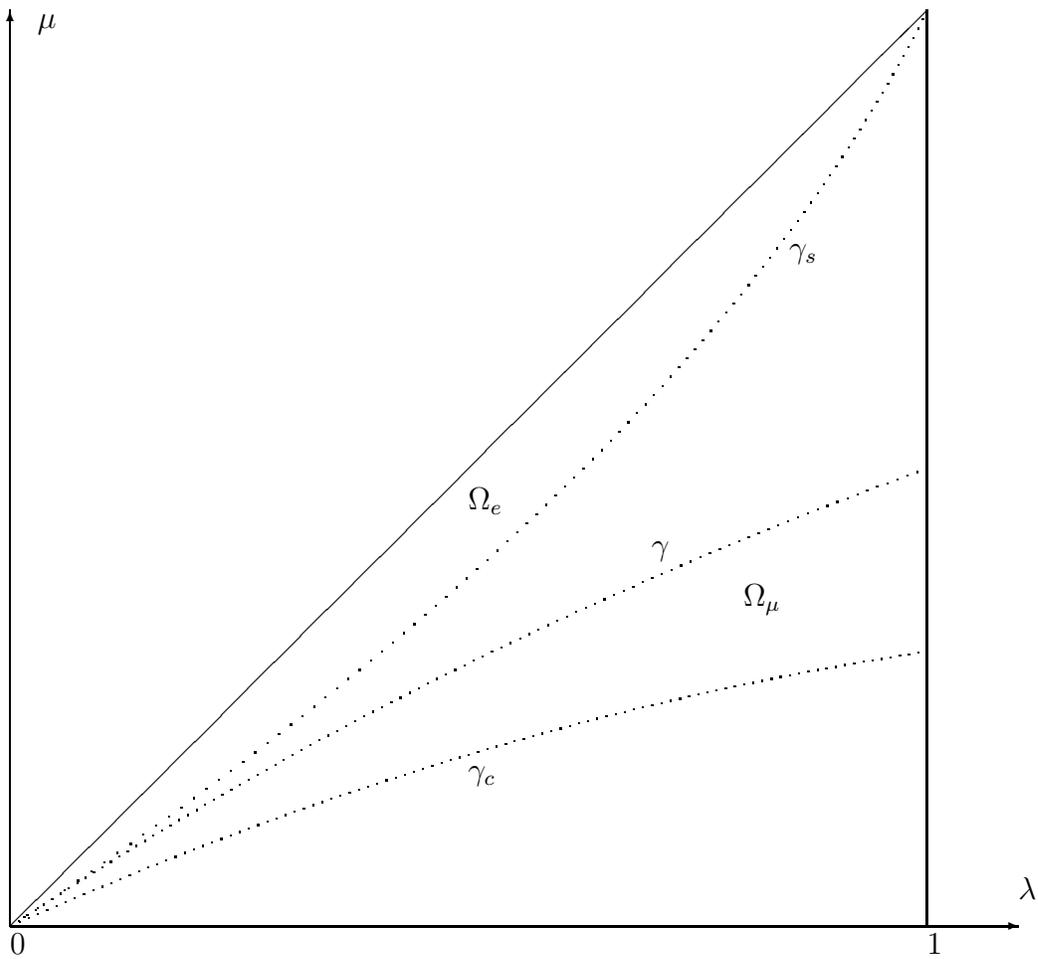
\begin{proposition}
\label{twocycles}
The vector spaces ${\cal A}_\mu$, ${\cal A}_e$, ${\cal A}_1$ of Abelian
integrals $\alpha_0I_0+\alpha_1I_1$ corresponding to the continuous
families of ovals ${\cal O}_\mu$ for $(\lambda,\mu)\in\Omega_\mu$,
${\cal O}_e$ for $(\lambda,\mu)\in\Omega_e$, and
${\cal O}_1$ for $\lambda\in\Omega_1\cup\Omega_1^*$, are not Chebyshev
in the respective intervals $\Sigma$.
\end{proposition}
{\bf Proof.} The assertion follows from local analysis.
Consider first ${\cal O}_\mu$. Then $\Sigma=(h_\mu, 0)$ and
$F(h_\mu)=\mu>0$, $F(0)=0$. By the proof of Proposition \ref{cyclicity},
we obtain that $F'(h_\mu)>0$ is equivalent to $a_1<0$. This coefficient
is easily calculated from the normal form (\ref{normal}) and the
result is
$$a_1=\frac{P'''(\mu)}{3P''(\mu)}=
\frac{6\mu^2-4\lambda\mu-4\mu+2\lambda}{3\mu(\lambda-\mu)(1-\mu)}$$
which proves the claim. The proof in the case of ${\cal O}_1$
is similar. The assertion for ${\cal O}_e$ follows from the
definition of $\gamma_s$. $\Box$

\vspace{1ex}
Making use of the whole information about the local behavior of $F$ near
the endpoints of $\Sigma$, one can formulate the following hypothesis.

\vspace{2ex}
\noindent
{\bf Conjecture 2.} (i).
The vector spaces ${\cal A}_\mu$, ${\cal A}_e$ are Chebyshev with accuracy
one in $\Sigma$ for the cases considered in Proposition \ref{twocycles},
and they are Chebyshev in $\Sigma$ for all the remaining cases.

(ii). The vector space ${\cal A}_1$ is Chebyshev with accuracy one in
$\Sigma$ for $\lambda\in\Omega_1$, it is Chebyshev with accuracy two in
$\Sigma$ for $\lambda\in \Omega_1^*$ and it is Chebyshev in $\Sigma$ for
the remaining cases.

(iii). Given $\alpha_0$ and $\alpha_1$, the total number of zeros of
the function $\alpha_0I_0+\alpha_1I_1$ relative to all the period annuli
is three and this bound is attained in $\Omega_e\cup\Omega_1^*$.

\vspace{2ex}
\noindent
It is plausible that Conjecture 2 is a consequence of Conjecture 1, except in
the case when the polynomial $P(x)$ has complex critical points.

The next problem in difficulty should be to study the space of Abelian
integrals of the first kind related to genus three hyperelliptic curves. We
hope that this problem can be settled by combining the methods of the present
paper with deformation arguments as in \cite{Gav4}. Properly understood,
this might lead to the solution of the problem for arbitrary genus $g$.

\vspace{2ex}
\noindent
{\bf Acknowledgments.} This research is partially supported by the
Ministry of Science and Education of Bulgaria under Grant no. 1003/2000.


\begin{thebibliography}{99}
\bibitem{cam}
%
N. A'Campo,
Le groupe de monodromie du d\'eploiement des singularit\'es isol\'ees de
courbes planes. I. {\it Math. Ann.} {\bf 213} (1975), 1--32.
%
\bibitem{arbarello} E. Arbarello, M. Cornalba, P.A. Griffiths, J. Harris,
{\it Geometry of Algebraic Curves}, vol. 1, Springer, New York, 1985.
%
\bibitem{Arnold} V.I. Arnold, S.M. Gusein-Zade, A.N. Varchenko,
{\it Singularities of Differentiable Maps}, vols. 1 and 2,
Monographs in Mathematics, Birkh\"auser, Boston, 1985 and 1988.
%
\bibitem{Arnold0} V.I. Arnold, Yu.S. Il'yashenko, Ordinary Differential
Equations, in: {\it Dynamical Systems} I, Encyclopaedia of Math. Sci.,
vol. 1, Springer, Berlin, 1988.
%
\bibitem{Arnold1} V.I. Arnold, {\it Geometrical Methods in the Theory of
Ordinary Differential Equations}, Springer, New York, 1988.
%
\bibitem{Arnold2} V.I. Arnold, Sur quelques probl\`emes de la
th\'eorie des syst\`emes dynamiques, {\it Topological Methods in Nonlinear
Analysis}, {\bf  4} (1994), 209--225.
%
\bibitem{Arnold1+} V.I. Arnold, Some unsolved problems in the theory of
differential equations and mathematical physics,
{\it Russian Math. Surveys} {\bf 44} (1989), no. 4, 157--171.
%
\bibitem{ten} V.I. Arnold, Ten problems, in: {\it Theory of
singularities and its applications}, {\it Adv. Soviet Math.} {\bf 1},
pp. 1--8, Amer. Math. Soc., Providence, 1990.
%
\bibitem{bri} E. Brieskorn, Die Monodromie der isolierten Singularit\"{a}ten
von Hyperfl\"{a}chen, {\it  Manuscripta Math.} {\bf 2} (1970), 103--161.
%
\bibitem{dim} A. Dimca, {\it Singularities and Topology of Hypersurfaces},
Springer, Berlin, 1992.
%
\bibitem{gas-li} A. Gasull, W. Li, J. Llibre, Zh. Zhang, Chebyshev property
of complete elliptic integrals and its application to Abelian integrals,
{\it Pacific J. Math.}  {\bf 202} (2002), no. 2, 341--362.
%
\bibitem{gav97} L. Gavrilov, Isochronicity of plane polynomial
Hamiltonian systems, {\it Nonlinearity} {\bf 10} (1997), no. 2, 433--448.
%
\bibitem{gav98a} L. Gavrilov, Petrov modules and zeros of Abelian
integrals, {\it Bull. Sci. Math.} {\bf 122} (1998), no. 8, 571--584.
%
\bibitem{gav99} L. Gavrilov, Abelian integrals related to Morse polynomials
and perturbations of plane Hamiltonian vector fields, {\it Ann. Inst. Fourier}
{\bf 49} (1999), no. 2, 611--652.
%
\bibitem{gav98} L. Gavrilov, Nonoscillation of elliptic integrals related
to cubic polynomials with symmetry of order three, {\it Bull. London Math.
Soc.} {\bf 30} (1998), no. 2,  267--273.
%
\bibitem{Gav4} L. Gavrilov, The infinitesimal 16th Hilbert problem in the
quadratic case, {\it Invent. Math.} {\bf 143} (2001), 449--497.
%
\bibitem{giv} A.B. Givental, Sturm's theorem for hyperelliptic integrals,
{\it Leningrad Math. J.} {\bf 1} (1990), 1157--1163.
%
\bibitem{gh} P.A. Griffiths, J. Harris, {\it Principles of algebraic
geometry}, Pure and Appl. Math., John Wiley and Sons, New York, 1978.
%
\bibitem{hil} D. Hilbert, Mathematische probleme,
{\it Gesammelte Abhandlungen} III, Springer-Verlag, Berlin (1935),
pp. 403--479.
%
\bibitem{nov-yak} D. Novikov and S. Yakovenko, Tangential Hilbert problem
for perturbations of hyperelliptic Hamiltonian systems, {\it Electron. Res.
Announc. Amer. Math. Soc.} {\bf 5} (1999), 55--65 (electronic).
%
%
\bibitem{Petrov86} G.S. Petrov, Number of zeros of complete elliptic integrals,
{\it Funct.Anal.Appl.} {\bf 18} (1984), 72--73; {\bf 20} (1986), 37--40;
{\bf 21 } (1987), 87--88; {\bf 22} (1988), 37--40;
{\bf 23} (1989), 88--89; {\bf 24} (1990), 45--50.
%
\bibitem{pont} L.S. Pontryagin, On dynamic systems close to Hamiltonian
systems, {\it Zh. Eksp. Teor. Fiz.} {\bf 4} (1934), 234--238 (Russian).
%
\bibitem{Roussarie1} R. Roussarie, {\it Bifurcation of Planar Vector Fields
and Hilbert's sixteenth Problem}, Progress in Mathematics, vol. 164,
Birkh{\"a}user, Basel, 1998.
%
\end{thebibliography}
\end{document}